\documentclass[12pt]{article}
\usepackage{amssymb,epsfig}
\newtheorem{defin}{Definition}
\newtheorem{prop}{Proposition}
\newtheorem{nt}{Remark}
\newtheorem{Th}{Theorem}

\newfont{\sdbl}{msbm9}
\newfont{\dbl}{msbm10 at 12pt}

\newcommand{\proof}{{\bf Proof\ }}

\newcommand{\oo}{{\cal O}}

\newcommand{\g}{{\cal G}}

\newcommand{\End}{\mathop {\rm End}}
\newcommand{\Lim}{\mathop {\rm lim}}

\newcommand{\sda}{{\mbox{\sdbl A}}}
\newcommand{\da}{{\mbox{\dbl A}}}
\newcommand{\dz}{{\mbox{\dbl Z}}}

\newcommand{\sdz}{{\mbox{\sdbl Z}}}

\newcommand{\f}{{\cal F}}
\newcommand{\lto}{\longrightarrow}

\begin{document}
\author{Denis Osipov
\footnote{supported by DFG-Schwerpunkt ''Globale Methoden in der
Komplexen Geometrie'', by RFBR grant no. 05-01-00455 and grant no.
04-01-00702, by grant of Leading Scientific Schools no. 9969.2006.1,
by INTAS grant 05-100000-8118, and by grant of Russian Science
Support Foundation.} }

\title{Adeles on  $n$-dimensional schemes and categories $C_n$.}
\date{}

\maketitle

\abstract{We consider categories $C_n$ which are very close to the
iterated functor $\mathop{\lim}\limits_{\longleftrightarrow}$,
which was introduced by A.A.Beilinson in \cite{B1}. We prove that
an adelic space on $n$-dimensional Noetherian scheme  is an object
of $C_n$.}

\section{Introduction}

In this note we want to introduce by induction some class of
infinite-dimensional vector spaces and morphisms between them. These
spaces depend on integer $n$, we call such spaces as $C_n$-spaces.

$C_0$-spaces are finite dimensional spaces. A first basic example
of $C_1$-space is the field of Laurent series $k((t))$. If $\oo(n)
= t^n k[[t]]$, then we have a filtration
$$
\ldots  \oo(m) \subset \oo(m-1) \subset \oo(m-2) \ldots
$$
and every factor space $\oo(m - k)/ \oo(m) $
is  a finite dimensional vector space over $k$.

We can consider the morphisms between two such spaces as continuous
linear maps. We remark that these morphisms can be described only in
terms of filtration $\oo_n$, without considering topology on
$k((t))$. From this point of view,  the space of adeles on an
algebraic curve has a structure of $C_1$-space, which is filtered by
partially ordered set of coherent sheaves on the curve.

We construct an iterated version of $C_1$-spaces, which we call a
$C_n$-space. But in our construction we do not consider the
structure of completion. So, we consider the filtered vector
spaces. For example, the discrete valuation field is also a
$C_1$-space, the space of rational adeles from~\cite{S} is also a
$C_1$-space.

The  constructions of similar categories were introduced also
in~\cite{B1}, \cite{Ka}, \cite{K}. Our construction of $C_n$ is very
close to the iterated functor
$\mathop{\lim}\limits_{\longleftrightarrow}$, which was introduced
by A.A.Beilinson in appendix to~\cite{B1}. The main difference is
that we consider noncompleted version of
$\mathop{\lim}\limits_{\longleftrightarrow}$, i.e., filtered spaces,
but with morphisms which come from
$\mathop{\lim}\limits_{\longleftrightarrow}$. From this point of
view the categories $C_n$ are rather closed to the dir-inv modules
which were considered by A. Yekutieli in~\cite{Y} for $n=1$.

The main result of this note is the theorem \ref{th1}, where we
prove that the space of Parshin-Beilinson adeles on an
$n$-dimensional Noetherian scheme $V$ over $k$ is a $C_n$-space,
which is filtered by partially ordered set of coherent sheaves on
the scheme. We calculate also the endomorphism algebra of
n-dimensional local field.

Adeles on algebraic surfaces were introduced by A.N.Parshin
in~\cite{P2}. A.A.~Beilinson generalized it to  arbitrary
Noetherian schemes in~\cite{B}. Adeles on higher-dimensional
schemes were applied to a lot of problems of algebraic geometry,
see~\cite{PF}.

A.N. Parshin has pointed out to me that the categories $C_n$ can be
useful for constructing of harmonic analysis on higher-dimensional
schemes and higher-dimensional adeles and local fields, see
\cite{P1}.

In this note we consider the categories $C_n$ and the schemes over a
field $k$. But all the constructions, for example, can be moved to
arithmetical schemes, where $C_0$ are finite abelian groups, $C_1$
are filtered abelian groups with finite abelian group factors and so
on.

\section{Categories $C_n$}
\subsection{Constructions}

\subsubsection{Objects in $C_n$}

\begin{defin} \label{def1}
We say that $(I, F, V)$
is a filtered $k$-vector space, if
\begin{enumerate}
\item $V$ is a vector space over the field $k$, \item $I$ is a
partially ordered set, such that for any $i,j \in  I$ there are
$k,l \in I$ with $k \le i \le l$ and $k \le j \le l$, \item $F$ is
a function from $I$ to the set of $k$-vector subspaces of $V$ such
that if  $i \le j $ are any from $I$, then $F(i) \subset F(j)$,
\item $ \bigcap\limits_{i \in I}  F(i) = 0$  and $
\bigcup\limits_{i \in I} F(i) = V $.
\end{enumerate}
\end{defin}

\begin{defin}
We say that a filtered vector space $(I_1, F_1, V)$ dominates another filtered vector space $(I_2, F_2, V)$
when there is a preserving order function $\phi: I_2 \to I_1$ such that
\begin{enumerate}
\item for any $i \in I_2$ we have $F_1(\phi (i)) = F_2 (i)$
\item for any $j \in I_1$ there are $i_1, i_2 \in I_2$ such that $ \phi(i_1) \le j \le \phi(i_2)$.
\end{enumerate}
\end{defin}

Now we define by induction the category of $C_n$-spaces and morphisms between them.

\begin{defin}
\begin{enumerate}
 \item The category  $C_0$ is  the category of finite-dimensional vector spaces over $k$
 with morphisms coming from  $k$-linear maps between vector spaces.
 \item The triple from $C_0$
 $$0 \lto V_0 \lto V_1  \lto V_2 \lto  0 $$ is admissible when it is an exact triple of vector spaces \mbox{.}
 \end{enumerate}
\end{defin}

Now we define the objects of the  category $C_n$ by induction.
We suppose that we have already defined the objects of the category $C_{n-1}$
and the notion of admissible triple in $C_{n-1}$.

\begin{defin} \label{def4}
\begin{enumerate}
\item Objects of the category $C_n$, i.e. $Ob(C_n)$, are filtered $k$-vector spaces $(I, F, V)$ with the following additional
structures
\begin{enumerate}
\item for any $i \le j \in I$  on the $k$-vector space $F(j) / F(i)$ it is given a structure $E_{i,j} \in Ob(C_{n-1})$,
\item for any $i \le j \le k \in I$
$$
0 \lto E_{i,j}  \lto E_{i,k}  \lto  E_{j,k} \lto 0
$$
is an admissible triple from $C_{n-1}$.
\end{enumerate}
\item Let $E_1 = (I_1, F_1, V_1)$, $E_2 = (I_2, F_2, V_2)$ and $E_3=(I_3, F_3, V_3)$ be from  $Ob(C_n)$. Then
we say that
$$
0 \lto E_1 \lto E_2 \lto E_3 \lto 0
$$
is an admissible triple from $C_n$ when
 the following conditions are satisfied
\begin{enumerate}
\item
$$
0 \lto V_1 \lto V_2 \lto V_3 \lto 0
$$
is an exact triple of $k$-vector spaces
\item \label{itaa}
the filtration $(I_1, F_1, V_1)$ dominates the filtration $(I_2, F'_1, V_1)$,
where $F'_1 (i) = F_2(i) \cap V_1$  for any $i \in I_2$,
\item \label{itbb}
the filtration $(I_3, F_3, V_3)$ dominates the filtration $(I_2, F'_3, V_3)$,
where $F'_3(i) =  F_2(i) / F_2(i) \cap V_1$,
\item for any $i \le j \in I_2$
\begin{equation} \label{trojkaa1}
0 \lto \frac{F'_1(j)}{F'_1(i)} \lto \frac{F_2(j)}{F_2(i)}  \lto \frac{F'_3(j)}{F'_3(i)} \lto 0
\end{equation}
is an admissible triple from $C_{n-1}$. (By definition of $Ob(C_n)$, on every vector space from triple~(\ref{trojkaa1})
it is given the structure of $Ob(C_{n-1})$).
\end{enumerate}
\end{enumerate}
\end{defin}

\subsubsection{Morphisms in $C_n$}

 By induction, we define now  the morphisms in the category $C_n$. We suppose that we have already defined  the morphisms
in $C_{n-1}$.

\begin{defin} \label{d1}
Let $E_1 = (I_1, F_1, V_1)$ and $E_2 = (I_2, F_2, V_2)$ be from $Ob(C_n)$.
Then $Mor_{C_n}(E_1, E_2)$ consists of elements $A \in Hom_k (V_1, V_2)$ such that the following conditions hold.
\begin{enumerate}
\item \label{i1} for any $i \in I_1$ there is an $j \in I_2$ such that $A (F_1(i)) \subset F_2(j)$,
\item \label{i2}  for any $j \in I_2$ there is an $i \in I_1$ such that $A (F_1(i)) \subset F_2(j)$,
\item \label{i3}  for any $i_1 \le i_2 \in I_1$ and $j_1 \le j_2 \in I_2$ such that $A (F_1(i_1)) \subset F_2(j_1)$
and $A (F_1(i_2)) \subset F_2(j_2)$ we have that the induced $k$-linear map
$$   \bar{A} : \frac{F_1(i_2)}{F_1(i_1)} \lto \frac{F_2(j_2)}{F_2(j_1)}
$$
is an element from
$$Mor_{C_{n-1}}(\frac{F_1(i_2)}{F_1(i_1)}, \frac{F_2(j_2)}{F_2(j_1)})$$
\end{enumerate}
\end{defin}

Now we want to prove that the compositions of so defined morphisms in $C_n$ will be again a morphism.
We need the following definition.

\begin{defin}
 Let $E_1, E_2$ be from $Ob(C_n)$.
\begin{enumerate}
\item
A $k$-linear map $C: E_1 \to E_2$ is an admissible
$C_n$-monomorphism
when it is the part of an admissible triple from $C_n$
$$
0 \lto E_1 \stackrel{C}{\lto} E_2 \lto E_3 \lto 0
$$
\item
A $k$-linear map $D: E_1 \to E_2$ is an admissible
$C_n$-epimorphism
when it is the part of an admissible triple
$$
0 \lto E_3 \lto E_1 \stackrel{D}{\lto} E_2 \lto 0
$$
\end{enumerate}
\end{defin}

We have the following proposition.
\begin{prop}
Let $E_1= (I_1, F_1, V_1)$, $E_2= (I_2, F_2, V_2)$, $E_1'$, $E_2'$ be from
$Ob(C_n)$ and  $A$ is a $k$-linear map form  $Hom_k(V_1, V_2)$.
\begin{enumerate}
\item \label{item1}
If $B: E_3 \to E_1$ is an admissible $C_n$-epimorphism then $A \in Mor_{C_n}(E_1, E_2)$ if and only if
$ A \circ B \in Mor_{C_n}(E_3, E_2)$.
\item \label{item2}
If $B: E_2 \to E_3$ is an admissible $C_n$-monomorphism, then $A \in Mor_{C_n}(E_1, E_2)$ if and only if
$ B \circ A \in Mor_{C_n}(E_1, E_3)$.
\item If the filtered vector space $E_1$ dominates the filtered vector space $E'_1$
and the filtered vector space $E_2$ dominates the filtered vector space $E'_2 $, then $A \in Mor_{C_n}(E_1, E_2)$ if and only if
$A \in Mor_{C_n}(E'_1, E'_2)$.
\item $Mor_{C_n}(E_1, E_2)$ is a $k$-linear subspace of $Hom_k(V_1, V_2)$.
\item If $E_3$ is an object of $C_n$, then
$$     Mor_{C_n}(E_2, E_3)   \circ   Mor_{C_n}(E_1, E_2)  \subset  Mor_{C_n}(E_1, E_3) $$
\end{enumerate}
\end{prop}
\proof.
The first two statements follow by induction on $n$.

The third statement follows from the first and the second statement.

The other statements follow by induction on $n$ using the previous statements.

We give the proof of the fifth statement.
Let $A \in Mor_{C_n}(E_1, E_2)$ and $B \in Mor_{C_n}(E_2, E_3)$. We have to prove that
$B \circ A \in Mor_{C_n}(E_1, E_3)$. We have to check for $B \circ A$ the conditions~\ref{i1}, \ref{i2}, \ref{i3}
of definition~\ref{d1}. Let $E_3 = (I_3, F_3, V_3)$.

For any $i_1 \in I_1$ there is $i'_2 \in I_2$ such that $A (F_1(i_1)) \subset F_2(i'_2)$.
For $i'_2 \in I_2$ there is $i'_3 \in I_3$ such that $B (F_2(i'_2)) \subset F_3 (i'_3)$.
Therefore $B \circ A (F_1 (i_1)) \subset F_3 (i'_3)$.

Analogously for any $j_3 \in I_3$ we find $j'_2 \in I_2$ such that $B (F_2(j'_2)) \subset F(j_3)$.
For $j'_2 \in I_2$ we find $j'_1 \in I_1$ such that $A (F_1(j'_1)) \subset F_2 (j'_2)$.
Then $B \circ A  (F_1(j'_1)) \subset F_3(j_3)$.

Now let $i_1 \ge j_1 \in I_1$ and $i_3 \ge j_3 \in I_3$ such that
$$ B \circ A (F_1(j_1)) \subset F_3(j_3)   \qquad \mbox{and}  \qquad B \circ A (F_1(i_1)) \subset F_3(i_3) \mbox{.}$$

Now we fix any $i^{''}_3 \in I_3$ such that $i^{''}_3 \ge i'_3$ and
$i^{''}_3 \ge i_3$. We fix any $j^{''}_1 \in I_1$ such that
$j^{''}_1 \le j_1$ and $j^{''}_1 \le j'_1$. Then, by
items~\ref{item1} and~\ref{item2} of this proposition, the map
induced by $B \circ A$ belongs to
$Mor_{C_{n-1}}(\frac{F_1(i_1)}{F_1(j_1)},
\frac{F_3(i_3)}{F_3(j_3)})$ if and only if the map induced by $B
\circ A$ belongs to $Mor_{C_{n-1}}(\frac{F_1(i_1)}{F_1(j^{''}_1)},
\frac{F_3(i^{''}_3)}{F_3(j_3)})$. But the last induced map is the
composition of the induced map by $A$ in
$Mor_{C_{n-1}}(\frac{F_1(i_1)}{F_1(j^{''}_1)},
\frac{F_2(i'_2)}{F_2(j'_2)})$ and of the induced map by $B$ in
$Mor_{C_{n-1}}(\frac{F_2(i_2)}{F_2(j'_2)},
\frac{F_3(i^{''}_3)}{F_3(j_3)})$. By induction, the composition of
morphisms from $C_{n-1}$ is a morphism from $C_{n-1}$. The
proposition is proved.

\subsection{Examples}

\subsubsection{Linearly locally compact spaces.} \label{lcs}
\begin{defin}
 A topological $k$-vector space $V$ (over a discrete field $k$) is linearly compact
 (see \cite[ch.2 ,\S 6]{L}, \cite[ch.III, \S 2, ex.15-21]{Bo})
when the following conditions hold
\begin{enumerate}
\item $V$ is complete and Hausdorff, \item $V$ has a base of
neighbourhoods
 of $0$ consisting of vector subspaces, \item each
open subspace of $V$ has finite codimension.
\end{enumerate}
\end{defin}

\begin{defin}
A topological $k$-vector space $W$ (over a discrete field $k$) is
linearly locally compact (see \cite{L}) when it has a basis of
neighborhoods of $0$ formed by linearly compact open subspaces.
\end{defin}

Any topological  linearly locally compact space is a $C_1$-space,
where filtration is given by linearly compact open subspaces.

Any field of discrete valuation is $C_1$-space. And the completion
functor gives us the linearly locally compact vector space.
 Moreover, if $E = (I, F, V)$ is a $C_1$-space, we take
$$
 \Phi_1(E) = \mathop{\Lim_{\rightarrow}}_{i \in I} \mathop{\Lim_{\leftarrow}}_{j \le i} F(i) / F(j) \mbox{.}
$$
The  space $\Phi_1 (E)$ is a  linearly locally compact space.
And all the linearly locally compact spaces can be obtained in a such way.

Moreover, we  define by induction on $n$ the functor of completion $\Phi_n$ from $C_n$-spaces to $C_n$-spaces.
We put
$$\Phi_n (E) =
\mathop{\Lim_{\rightarrow}}_{i \in I} \mathop{\Lim_{\leftarrow}}_{j \le i}  \Phi_{n-1} ( F(i) / F(j))
$$
where $E= (I, F, V) $ is a $C_n$-space. From properties of $C_n$ we obtain that this functor $\Phi_n$ is well-defined.

\subsubsection{Contragredient space}
We will define by induction on $n$ a contravariant functor of
contragredient space $D_n$ from $C_n$-spaces to $C_n$-spaces. Let $E
=(I, F, V)$ be a $C_n$-space.

If $n=0$, then $D_0 (V) = V^*$, where $V^*$ is the dual space to
$V$.

If $n \ge 1$, then we put
$$
D_n (V) = \mathop{\Lim_{\rightarrow}}_{j \in I}
\mathop{\Lim_{\leftarrow}}_{i \ge j}  D_{n-1} ( F(i) / F(j))
\mbox{.}
$$

Then $D_n(E) = (I^0, F^0, D_n(V) )$ is a $C_n$-space, where $I^0$ is
a partially ordered set, which has the same set as $I$, but with the
inverse order then $I$, and
$$F^0 (j) =
\mathop{\Lim_{\leftarrow}}_{i \le j \in I^0} D_{n-1} ( F(i) / F(j))
\mbox{.}
$$

It is easy to see by induction on $n$ that $D_n(V) \subset V^*$, the
functor $D_n$ maps admissible triples to admissible triples, and
$D_n(D_n (E)) = \Phi_n(E)$.

\subsubsection{Adelic space.}

The next  example of $C_n$-space is the space of adeles on a
Noetherian $n$-dimensional scheme $V$ (see~\cite{B}, \cite{H},
\cite{Osi}).

Let $P(V)$ be the set of points of the scheme $V$.
Consider $\eta, \nu \in P(V)$. Define $\eta \ge \nu$
if $\nu \in \bar{\{ \eta\}}$. $\ge$ is a half ordering on $P(X)$.
Let $S(V)$ be the simplicial set induced by $(P(X), \ge)$,
i.e.
$$S(V)_m = \{ (\nu_0, \ldots, \nu_m) \mid \nu_i \in P(V); \nu_i \ge \nu_{i+1} \}$$
is the set of $m$-simplices of $S(V)$
with the usual boundary and degeneracy maps.

In \cite{B}, \cite{H} for any $K \subset S(V)_m$ it was
constructed the  space $\da(K, \f)$ for any quasicoherent sheaf
$\f$ on $V$ such that
$$
\da(K, \f) \subset
\prod_{\delta \in K} \da(\delta, \f)
$$

\begin{Th} \label{th1}
Let $V$ be a Noetherian $n$-dimensional scheme over $k$, $K
\subset S(V)_n$, and $\f$ be a coherent sheaf on $V$. Then the
adelic space $\da(K, \f)$  has a structure of $C_n$-space.
\end{Th}
\proof.
We denote
 $$K_0 =  \{\eta \in S(V)_0  \mid   (\eta > \eta_1 \ldots > \eta_n ) \in K  \quad
\mbox{for some} \quad \eta_i \in P(V) \} \mbox{.} $$

For $\eta \in K_0$ we denote
$$
{}_{\eta} K = \{ (\eta_1 > \ldots > \eta_{n}) \in S(V)_{n-1}
 \mid (\eta > \eta_1 \ldots  > \eta_n) \in K \} \mbox{.}
$$

We have
$$
\da(K,\f) = \prod_{\eta \in  K_0} \da({}_{\eta} K ,\f_{\eta}) \mbox{.}
$$

 If $(I_1, F_1, V_1)$ and $(I_2, F_2, V_2)$ are $C_n$-spaces,
then $(I_1 \times I_2, F_1 \times F_2, V_1 \times V_2)$ is  a $C_n$-space as well.
Moreover, any finite product of $C_n$-spaces is a $C_n$-space in the same way.
The set $K_0$ is  finite,
therefore it is enough to define a $C_n$-structure on  $\da({}_{\eta} K ,\f_{\eta})$ for every $\eta \in K_0$.

For $\eta \in K_0$ we define a partially ordered set
$$ I_{\eta}(\f) = \{ \g \subset \f_{\eta} \mid  \quad  \g \quad \mbox{is a coherent sheaf on }
\bar{\eta}, \quad \g_{\eta} = \f_{\eta} \}  \mbox{,}$$
which is ordered by inclusions of sheafs.

The functor $\da({}_{\eta} K, \;)$ is an exact functor. Therefore for any $\g \in I_{\eta}(\f)$ we have an embedding
$$\da({}_{\eta} K,\g) \lto   \da({}_{\eta} K,\f_{\eta}) \mbox{.}$$

If $\g_1 \subset \g_2$ are from $I_{\eta} (\f)$,
then from the exactness of the functor $\da({}_{\eta} K, \;)$ we have
$$
\da({}_{\eta} K, \g_2) /
\da({}_{\eta} K,  \g_1) =
\da({}_{\eta} K, \g_2 / \g_1)
\mbox{.}
$$

We have $(\g_2 / \g_1)_{\eta} = 0$.
Therefore $\g_2 / \g_1 $
is a coherent sheaf on some subscheme $Y$ of dimension $n-1$.
Therefore
$\da({}_{\eta} K, \g_2 / \g_1) =
\da ( S(Y) \cap {}_{\eta} K , \g_2 / \g_1 )
$.
We apply induction on dimension of scheme to the $n-1$-dimensional scheme $Y$.
Therefore we have that $\da ( S(Y) \cap {}_{\eta} K, \g_2 / \g_1 )$
is a $C_{n-1}$-space.

We check now that this $C_n$-structure is well defined.
It is enough to prove only
 conditions~\ref{itaa} and \ref{itbb} of definition~\ref{def4}.
 They follow from the following three statements.

  The first statement  claims
 that every exact triple of quasicoherent sheaves on $V$
 is a direct limit of exact triples of coherent sheaves on $V$, see~\cite[lemma 1.2.2]{H}.

 The second statement claims that on any irreducible Noetherian scheme $X$ for any exact triple
 of coherent sheaves
 $$
 0 \lto  \f_1 \lto \f_2 \lto \f_3 \lto 0
 $$
 and any coherent subsheaves $\g_1 \subset \f_1$ and $\g_3 \subset \f_3$ such that
 $(\f_1 / \g_1)_{\eta} = 0$ and $(\f_3 / \g_3 )_{\eta}= 0$, where $\eta $
 is the general point of $V$, there exists  a coherent subsheaf $\g_2 \subset \f_2$
 such that $ (\g_2 \cap \f_1) \subset \g_1 $, $ (\g_2 / (\g_2 \cap \f_1)) \subset \g_3 $
 and $(\f_2 / \g_2)_{\eta} = 0$. It is enough to construct $\f_2$ locally, where this sheaf exists by
 the Artin-Rees lemma.

 The third statement claims that for any two quasicoherent subsheaves
 $\f_1$, $\f_2$ of a quasicoherent sheaf $\f_3$ on $V$ we have
 $$
 \da (K, \f_1 \cap \f_2) = \da (K, \f_1) \cap \da (K, \f_2) \mbox{,}
 $$
 $$
 \da (K, \f_1 + \f_2) = \da (K, \f_1) + \da(K, \f_2) \mbox{,}
 $$
 where the intersection and sum is taken inside of $\da (K, \f_3)$.
 It follows from exactness of the functor $\da (K, \,)$ and the following commutative
 diagram:
 $$
\begin{array}{ccc}
\frac{\sda (K, \: \f_1 + \f_2)}{\sda (K, \: \f_1)} &  =  &
\frac{\sda(K, \: \f_2)}{\sda (K, \, \f_1 \cap \f_2)} \\
\uparrow & & \downarrow \\
\frac{\sda (K, \: \f_1) + \sda (K, \: \f_2)}{\sda (K, \: \f_1)} & =
& \frac{\sda (K, \: \f_2)}{\sda (K, \: \f_1) \cap \sda (K, \: \f_2)}
\quad \mbox{,}
\end{array}
 $$
  which gives that the vertical arrows are isomorphisms.

 The theorem is proved.

\begin{nt}  {\em
For the smooth surface $V$ a $C_2$-structure on $\da(K,\f)$ can be
defined by filtration of Cartier divisors. (We use that the
Cartier divisors coincide with the Weil divisors for the smooth
varieties). }
\end{nt}

\begin{nt} {\em
The structure of $C_n$ space can be defined on introduced
in~\cite[\S 5.2]{H}  spaces of rational adeles $a(K, \f)$ as well.
Then the functor of completion $\Phi_n$ (see section~\ref{lcs})
applied to the space of rational adeles gives the space  $\da(K,
\f)$.}
\end{nt}

\subsection{The endomorphism algebra of $n$-dimensional local field}
We suppose now that $\delta = (\eta_0 > \ldots > \eta_n) \in S(V)_n$
on the $n$-dimensional scheme $V$ has the following property:
$\eta_n$ is a smooth point on every scheme $\bar{\eta_i}$.

Then $ (\oo_V)_{\delta} =
\da(\delta, \oo_V)=
 k'((t_n)) \ldots ((t_1))$ is an $n$-dimensional local field,
where $k' = k(\eta_n)$.
We demand that  the local parameters $t_i \in \widehat{(\oo_V)}_{\eta_n}$ for any $i$.

Then we  define the filtration on $k' ((t_n)) \ldots ((t_1))$ by $ E_l = t_1^{l} k'((t_n)) \ldots ((t_2)) [[t_1]] $ for
$l \in \dz$. On each factor $E_{l_1} / E_{l_2}$ of this filtration  we define the new filtration
given in $E_{l_1} / E_{l_2}$ by images  of $ t_2^{m} t_1^{l_1} k'((t_n)) \ldots ((t_{3})) [[t_2, t_1]]$ for $m \in \dz$ and so on.

We obtained
 the
structure of $C_n$-space on $(\oo_V)_{\delta}$. And the structure of $C_n$-space constructed on $(\oo_V)_{\delta}$ in theorem~\ref{th1}
 dominates the constructed now structure of $C_n$-space.

Now let $K = k((t_n)) \ldots ((t_1))$. We define the $k$-algebra
 $$End_K = Mor_{C_n} (K, K)$$
Let $\bar{K} = k((t_{n})) \ldots ((t_2))$. Then $K = \bar{K}((t_1))$.
For any element $A \in \End_K$ we consider the matrix $\{(A_{ij})_{i,j \in \sdz} \mid A_{ij} \in End_{k}(\bar{K}) \}$
given by
$$
A (x t^i_n) = \sum_j A_{ij} ( x) t^j_n  \quad \mbox{with} \quad x
\in \bar{K} \mbox{.}
$$
\begin{prop}
An endomorphism $A \in End_k(K)$ belongs to $End_K$ if and only if the following conditions are satisfied.
\begin{enumerate}
\item
There is a nondecreasing function $a: \dz \to \dz$, $a(i) \to \infty$ when $i \to \infty$
such that for $j < a(i)$ all elements $A_{ij} = 0$.
\item Any element $A_{ij}$ belongs to $End_{\bar{K}}$.
\end{enumerate}
\end{prop}
\proof follows directly from the definition of morphisms between $C_n$-spaces and the definition of
$C_n$-structure on $K$. The proposition is proved.

We consider on $K$ the topology of $n$-dimensional local field (see \cite{P}).
Let $End^{c}_k(K)$ be the algebra of continuous $k$-linear endomorphisms of $K$.
\begin{prop}
\begin{enumerate}
\item \label{iit1} $End_K \subset End^{c}_k(K) $
\item \label{iit2} If $n =1$, then  $End_K = End^{c}_k(K)$.
\item \label{iit3} If $n > 1$, then $End^{c}_k(K)$ bigger then $End_K$.
\end{enumerate}
\end{prop}
\proof. We denote by $\oo_K = \bar{K}[[t_1]]$. The base of
neighbourhoods  of $0$ in $\bar{K}((t_1))$ consists of the
following  $k$-vector subspaces:
$$
\sum_i U_i t_1^i +  \oo_K t_1^m  \mbox{,}
$$
where $U_i$ are open $k$-vector subspaces from $\bar{K}$
and $m$ is an integer.

Then statement~\ref{iit2} is the same as condition~\ref{i2} of definition~\ref{d1}.
Statement~\ref{iit1} follows by induction on $n$ from the definition~\ref{d1},
for $n=2$ see~\cite[lemma~2]{O}.
Now we give an example for statement~\ref{iit3}.

Let $n=2$. Then we consider a $k$-linear map
$\phi:  K \to K $ defined as following.
For $U = k[[t_2]] ((t_1)) + k((t_2))[[t_1]])$ we put $\phi (U) = 0$
and the induced  map $\phi : K/U \to K $ we put on monomials $t_2^l t_1^m$ by the rule:
if $m=-1 $ then $\phi(t_2^l t_1^m) = t_2^l  t_1^{m+l}$,
for other monomials we put $\phi = 0$.
Then $\phi$ is continuous, but $\phi$ is not from $End_K$.
The proposition is proved.

\begin{prop} \label{pro}
\begin{enumerate}
\item For any $m$ we have an embedding $End_K^{\oplus m}  \hookrightarrow End_K$.
\item For any $m$ we have an embedding $ gl(m, K) \hookrightarrow End_K$.
\end{enumerate}
\end{prop}
\proof
We have $K = \bar{K}((t_n))$.
 Let $e_i$,    $1 \le i \le m$ be the standard basis of $K^{\oplus m}$, i.e. $K^m = \oplus K e_m$.
 We consider a $\bar{K}$-isomorphism $\phi$ of $\bar{K}$-vector spaces
$$
\phi : K^{\oplus m} \to K \quad  \quad  \phi(t_n^j e_i) = t_n^{mj +i -1}  \qquad 1 \le i \le m \mbox{.}
$$
This isomorphism induces an isomorphism
$$Mor_{C_n} (K^{\oplus m}, K^{\oplus m})  \hookrightarrow End_K \mbox{.} $$
Now the proposition follows from the natural  embeddings:
$$End_K^{\oplus m}  \hookrightarrow Mor_{C_n} (K^{\oplus m}, K^{\oplus m}) \qquad
\mbox{and} \qquad gl(m, K) \hookrightarrow  Mor_{C_n} (K^{\oplus m}, K^{\oplus m}) \mbox{.}$$

\begin{nt} {\em
From proposition~\ref{pro} we have embeddings of
$$K^* \hookrightarrow End_K$$
and toroidal Lie algebras
$$gl (m, k[t_1, \ldots, t_n, t_1^{-1}, \ldots, t_n^{-1}])
\hookrightarrow   End_K  \mbox{.} $$
}
\end{nt}

\section*{Acknowledgments}
I am very grateful to A.N.Parshin  for interesting and valuable
comments and the stimulation the author to write this note.

I am also grateful to H.Kurke,  A.Zheglov and I.Zhukov  for
interesting comments. This note was  prepared during my stay at the
Humboldt University of Berlin. I am grateful to the Humboldt
University of Berlin  for the hospitality.

\noindent Steklov Mathematical Institute, \\
Gubkina str. 8, \\
119991, Moscow, Russia \\
e-mail  ${d}_{-} osipov@mi.ras.ru$

\end{document}